\documentstyle[12pt]{article}

\begin{document}

\begin{center}
{\large \bf On dimension reduction in the K\"ahler-Ricci flow }
\end{center}
\begin{center}
\begin{tabular}{ccc}
Huai-Dong Cao \\
\end{tabular}
\end{center}
\footnotetext{Research supported in part by NSF grant DMS-0206847. }
\begin{abstract}
We consider dimension reduction for solutions of the K\"ahler-Ricci flow with nonegative bisectional curvature. When the complex dimension $n=2$, we prove an optimal dimension reduction theorem for complete translating K\"ahler-Ricci solitons with nonnegative bisectional curvature. We also prove a general dimension reduction theorem for complete ancient solutions of the K\"ahler-Ricci flow with nonnegative bisectional curvature on noncompact complex manifolds under a finiteness assumption on the Chern number $c^n_1$.  
\end{abstract}

\vspace{1 cm}

\section{Introduction}

In minimal surface theory and harmonic map theory, there is a well-known 
general principle of dimension reduction which is very useful in studying singularities. As is well-know, the basic idea is that after having taken 
a limit of a sequence of dilations to model a singularity (blow-up), one 
should study this limit by next taking a sequence of origin going out to 
infinity and shrinking back down (blow-down) to get a new limit of lower dimension. In [11], Richard Hamilton remarkably showed us how to carry out 
this general idea of dimension reduction to solutions of the Ricci 
flow, which is a nonlinear parabolic system, on Riemannian
manifolds with nonnegative curvature operator and local
injectivity radius estimate (cf. Section 22 of [11]), and used it
to prove important results about singularity formations of the
Ricci flow. See also the recent work of Chen-Tang-Zhu [7] for a
very nice application of Hamilton's dimension reduction result to
the K\"ahler-Ricci flow on noncompact complex surfaces with
positive bisectional curvature and maximal volume growth. For
studying formations of singularities of the K\"ahler-Ricci flow in
complex dimension $n\geq 2$, it has been very desirable but a
challenge to extend Hamilton's dimension reduction for the Ricci flow to
solutions of the K\"ahler-Ricci flow with positive holomorphic
bisectional curvature (wich is a weaker curvature assumption). 
 In the proof of Hamilton's dimension reduction result for the Ricci flow, the assumption of nonnegative curvature operator is used in an essential way in both finding a flat curvature direction of the blow-down limit (in terms of 
Theorem 21.6 of [11] on the finiteness of remote curvature bumps),
and the splitting of the flat factor. In the K\"ahler case, the
splitting theorem for complete solutions of the K\"ahler-Ricci
flow with nonnegative bisectional curvature (Theorem 2.1) has been
known to the author for a number of years. The main difficulty has
been about how to produce a flat direction in blow-down limits to
allow the splitting when we only have nonnegative bisectional
curvature. In this paper we resolve this difficulty and show how
the dimension reduction of Hamilton for the Ricci flow can be
extended to the K\"ahler case for $2$-dimensional translating K\"ahler-Ricci 
solitons (Theorem 1.1). We also prove a general dimension reduction result for ancient solutions of the K\"ahler-Ricci flow with nonnegative bisectional curvature under a certain finiteness assumption of the Chern number integral $c^n_1$ (Theorem 1.2). When complex dimension $n=2$, we are able to show that the second Chern number
integral is well defined (Proposition 4.1) for Type II singularity
models (or equivalently translating K\"ahler-Ricci solitons) of nonnegative bisectional curvature. By exploring the relation between the integrands of the Chern numbers 
$c^2_1$ and $c_2$, we can prove Theorem 1.1 without the finiteness assumption on $c_1^2$.

A complete solution to the K\"ahler-Ricci flow
\begin{eqnarray}
{\partial \over {\partial t}}g_{i\bar j}(x,t)= -R_{i\bar j}(x,t)
\end{eqnarray}
is a time-dependent family of complete K\"ahler metrics $g=\sum g_{i \bar j}(x,t)dz^i dz^{\bar j}$
on a complex manifold $X^n$ (either compact or noncompact) and some time interval  satisfying Eq. (1), where $R_{i \bar j}(x,t)$ denotes the Ricci tensor of the metric $g$ at time $t$. \\

\noindent {\bf Definition} (cf. Section 16 of [11]) A solution $g$ to the
K\"ahler-Ricci flow (1), where either the complex manifold $X$ is
compact or at each time $t$ the metric $g$ is complete with bounded and
nonnegative bisectional curvature, is called a {\em singularity
model} or {\em limit solution} of Eq.(1) if it is non flat and of
one of the following three types:

\begin{enumerate}
\item[]  Type I: $X$ is either compact or noncompact and $g$ exists on an ancient time interval $-\infty<t<\Omega$ for some constant $\Omega$ with $0<\Omega<\+\infty$, and the scalar curvature $R$ satisfies the inequality $$R(x,t)\leq\Omega/(\Omega -t)$$ everywhere with equality at some origin point $O$ at time $t=0$.

\item[]  Type II : $X$ is noncompact and $g$ exists on the eternal time interval $-\infty<t<+\infty$, and  $$R(x,t)\leq 1$$ everywhere with equality at some origin point $O$ at time $t=0$.

\item[]  Type III: $X$ is noncompact and $g$ exists on the time interval $-\Omega<t<\infty$ for some constant $\Omega$ with $0<\Omega<\+\infty$, and  $$R(x,t)\leq \Omega/(\Omega +t)$$ everywhere with equality at some origin point $O$ at time $t=0$.
\end{enumerate}
 The three types of singularity models above typically
arise as limit of blow-ups of maximal solutions to Eq.(1) of the
corresponding type. We remark that
Type II limit solutions of nonnegative bisectional curvature are necessarily
translating K\"ahler-Ricci solitons by our previous work [3] (see also Prop. 3.1). Note also
Type III limit solutions have the nice property that their
curvature decays in $t$ linearly as $t$ approaches the infinity
which in turn implies almost quadratic decay of curvature in space
in some average sense (cf. [7]). This makes Type III limit
solutions very special. In this paper, we will mainly focus our
attention on Type I and Type II singularity models. Our first
result is an optimal dimension reduction theorem for $2$-dimensional 
Type II singularity models:\\

\noindent{\bf  Theorem 1.1}  {\em Let $g$ be a Type II singularity
model, or equivalently translating K\"ahler-Ricci soliton, of the K\"ahler-Ricci 
flow (1) on a noncompact $2$-dimensional complex surface $X^2$ and the eternal time
interval $-\infty <t<\infty$, complete with bounded and
nonnegative bisectional curvature at each time. If $g$ satisfies
the local injectivity radius estimate
\begin{eqnarray}
 inj_X (x,t)\geq \beta/\sqrt{R(x,t)}
\end{eqnarray}
for all $x\in X^2$ and some constant $\beta >0$ independent of
$x$, then there exists a sequences of dilations of $g$ which
converges to a limit $\widehat g$, which is again a complete
solution to the K\"ahler-Ricci flow (1) on some noncompact complex
surface $\widehat X^2$ and ancient time interval $-\infty <t<
\Omega $, with nonnegative bisectional curvature.  Moreover,
$(\widehat X^2, \widehat g)$ splits as a quotient of
the product ${\bf C}\times {\bf CP}^1 $ of flat complex plane ${\bf C}$ and
the complex projective plane ${\bf CP}^1$. }\\

Whenever a Riemannian manifold $(M^m, g)$ of real dimension $m$
has nonnegative Ricci curvature, the ratio $V(B_s)/s^{m}$, where
$V(B_s)$ is the volume of geodesic ball of radius $s$ centered at
some origin point $O\in M$, is monotone decreasing in $s$ by the volume comparison theorem. So one
can consider the asymptotic volume ratio (cf. Hamilton [11], P.74)
\[\nu_{M}=\lim_{s\rightarrow \infty} V(B_s)/s^{m}.\] The
definition of $\nu_{M}$ is independent of the choice of the origin
$O$. Note that $\nu_M>0$ if and only if $(M^m, g)$ has Euclidean
volume growth.

\medskip

\noindent{\bf  Corollary 1.1}  {\em There exists no
$2$-dimensional Type II singularity models with bounded and nonnegative bisectional curvature, and with Euclidean volume
growth. In other words, the asymptotic volume ratio $\nu$ of any $2$-dimensional 
Type II singularity model with nonegative bisectional curvature  must be zero.}\\

\medskip

\noindent{\bf Remark}: Corollary 1.1 is also implicitly proved
in [7] when $g$ is assumed to have nonnegative curvature operator.

Our second result is a general dimension reduction for ancient solutions:\\

\noindent{\bf Theorem 1.2}  {\em Let $g$ be a complete ancient solution
to the K\"ahler-Ricci flow on a noncompact complex manifold $X^n$
($n\ge 2$) and time interval $-\infty <t< T (0<T\leq \infty)$, with bounded and nonnegative 
bisectional curvature, and positive Ricci curvature at each time. Assume that the following
conditions are satisfied:

\begin{enumerate}
\item[(i)] Finite asymptotic scalar curvature ratio (ASCR):
$$A\equiv\limsup_{s\rightarrow \infty} Rs^2=\infty$$ at $t=0$. Here $s$ denotes the distance function on $X$ from some fixed origin point;

\item[(ii)] Local injectivity radius estimate: there exists a constant $\beta >0$ such that
$\forall x\in X^n$ and $t$, 
\[inj_X (x,t)\geq
\beta/\sqrt{R(x,t)};\]

\item[(iii)] Finite Chern number: $ c_1^{n}(X)\equiv\int_X
Rc^{n}<\infty$.

\end{enumerate}
Then there exists a sequences of dilations of $g$ which converges
to a limit $\widehat g $,  which is again a complete solution to
the K\"ahler-Ricci flow (1) on some noncompact complex manifold
$\widehat X^n$ and ancient time interval $-\infty <t< \Omega $, with
nonnegative bisectional curvature.  Moreover, the limit $(\widehat
X^n, \widehat g)$ splits as a quotient of
a product ${\bf C^k}\times \widehat N^{n-k}$ with $k\ge 1$ flat in the
direction of ${\bf C}^k$, and where the interesting factor $\widehat N^{n-k}$
has positive Ricci curvature, and is either a Type I limit solution, or Type II
limit solution with $ c_1^{n-k}(\widehat N)=\infty$. }\\

\noindent {\bf Remark}: We conjecture that the only possible
factor $\widehat N$ in Theorem 1.2 is of Type I and compact. In
fact, we conjecture that there is no noncompact Type I limit
solutions with nonnegative bisectional curvatue and positive Ricci
curvature. Note that Hamilton (Section 26 of [11]) proved that in
complex dimension one, the only Type I ancient solution is either
the complex projective plane ${\bf CP}^1$, or flat complex plane
${\bf C}$ and its quotients.

\medskip

\noindent{\bf Remark}: Note that if $g$ comes as a blow-up
limit of a solution to the K\"ahler-Ricci flow (1) on a compact
K\"ahler manifold then condition (iii) in Theorem 1.2 is automatically 
satisfied. Moreover, in this case the assumption on the local 
injectivity radius in both Theorem 1.1 and Theorem 1.2 also holds according to a
very recent result of Perelman (cf. Corollary 4.3 of [15]). In case $g$ 
is a Type II limit solution, very likely condition (iii) could be
removed (when $n=2$, this is indeed the case as seen in Theorem 1.1).

\vfill \eject

\medskip

\noindent{\bf Acknowledgment}: A large part of the paper was written up while the
author was visiting the National Center of Theoretical Sciences in Hsinchu, Taiwan 
in July, 2002 where he also attended the five-week long summer workshop on Geometric Evolution Equations. The author would like to thank the Center for the hospitality and support during his visit. He would also like to thank Ben Chow, Richard Hamilton, Luen-fai Tam, Lei Ni and Jiaping Wang for their interest in this work.

\section{A splitting theorem for the K\"ahler-Ricci flow}

In this section we state and prove a splitting theorem for solutions of the K\"ahler-Ricci flow
(1) with nonnegative bisectional curvature. The result is a natural analogue to the splitting
theorem of Hamilton [9] for solutions of the Ricci flow on Riemannian manifolds with
nonnegative curvature operator. This splitting theorem will be useful in subsequent sections,
in particular in the proof of our dimension reduction theorems. \\

\noindent{\bf Theorem 2.1}   {\em Let $g$ be a complete solution
of the K\"ahler-Ricci flow (1) on a noncompact simply connected
complex manifold $X^n$ of dimension $n$ and some open time interval $I$, with bounded and
nonnegative holomorphic bisectioanl curvature. Then either $g$ is of positive Ricci curvature
for all $x\in X$ and all $t\in I$, or $(X, g)$ splits holomorphically isometrically into a
product $C^k \times N^{n-k}$ ($k\geq 1$) flat in $C^k$ direction and $N$ being of nonnegative
holomorphic bisectioanl curvature and positive Ricci curvature.}\\

\noindent{\bf Proof}: The proof can proceed essentially along the same line as in [9]. We claim that for any $t_0$, there exists a time interval $t_0<t<t_0+\delta$, on which the rank of the Ricci tensor $R_{i\bar j}$ is constant and the null space of $R_{i\bar j}$ is invariant under parallel translation and invariant in time. First, recall that the Ricci tensor satisfies the evolution equation:
\begin{eqnarray}
{\partial \over {\partial t}}R_{i\bar j}= \Delta R_{i\bar j} + \Phi_{i\bar j},
\end{eqnarray}
where $\Phi_{i\bar j}=R_{i\bar j k\bar l}R_{l\bar k}$.
Let $0\leq \sigma_1\leq \sigma_2\leq \cdots\leq \sigma_n$ be the eigenvalues of the Ricci tensor. Then
$\sigma_1+\cdots +\sigma_k$ is a concave function of $R_{i\bar j}$ and is invariant under parallel translation, since
$$ \sigma_1+\cdots +\sigma_k =\inf \{tr(R_{i\bar j}|E: E\subset T_X \  \mbox{is any subspace of dim k}\}.$$
Note that dim of the null space of  $R_{i\bar j}$ is $\geq k$ if and only if $\sigma_1+\cdots +\sigma_k=0$. If
$\sigma_1+\cdots +\sigma_k>0$ at one point at $t=t_0$, then by the strong maximum principle, it is positive everywhere for all $t>t_0$. So it follows that rank of $R_{i\bar j}$ remains constant on some time interval $t_0<t<t_0+\delta$.
\eject\vfill
Next, let $v$ be any smooth section of the holomorphic tangent bundle $T_X$ in the null space of $R_{i\bar j}$
on $t_0<t<t_0+\delta$. Then
\[ 0={\partial \over {\partial t}}(R_{i\bar j}v^iv^{\bar j})=({\partial \over {\partial t}}R_{i\bar j})v^iv^{\bar j}
+R_{i \bar j}({\partial v^i\over {\partial t}}v^{\bar j} + v^i {\partial v^{\bar j}\over {\partial t}}).\]
Since $R_{i\bar j}v^i=R_{i\bar j}v^{\bar j}=0$, the last term vanishes. Also
\begin{eqnarray*}
0&=&\Delta (R_{i\bar j}v^iv^{\bar j})=(\Delta R_{i\bar j})v^iv^{\bar j} + D_kR_{i\bar j}(v^{\bar j}D_{\bar k}v^i+v^iD_{\bar k}v^{\bar j}) \\
&+&D_{\bar k}R_{i\bar j}(v^{\bar j}D_{k}v^i+v^iD_{k}v^{\bar j})+R_{i\bar j}(D_kv^iD_{\bar k}v^{\bar j} +D_{\bar k}v^i D_kv^{\bar j}) + R_{i\bar j}(v^{\bar j}\Delta v^i +v^i\Delta v^{\bar j})
\end{eqnarray*}
 and again the last term disappears. Since
\[0=D_k(R_{i\bar j}v^i)=(D_k R{i\bar j})v^i + R_{i\bar j}D_kv^i,\]
and
\[0=D_k(R_{i\bar j}v^{\bar j})=(D_k R{i\bar j})v^{\bar j} + R_{i\bar j}D_kv^{\bar j},\] etc.,
 we get from the evolution equation (3)
\[R_{i\bar j}(D_kv^iD_{\bar k}v^{\bar j}+D_{\bar k}v^iD_kv^{\bar j})+\Phi_{i\bar j}v^iv^{\bar j}=0.\]

Since $R_{i\bar j}\geq 0$ and $\Phi_{i\bar j}=R_{i\bar j k\bar l}R_{l\bar k}\geq 0$, we must have $v$ also in the null space of $\Phi_{i \bar j}$ and $D_kv^i$, $D_kv^{\bar j}$ etc in the null spaces of $R_{i\bar j}$
for all k. This shows that the null space of the Ricci tensor is invariant under parallel translation and
null($R_{i\bar j}$)$\subset$ null($\Phi_{i\bar j}$).

Finally, to see the null space of the Ricci tensor is also invariant in time, note first that $\Delta v^i$ and $\Delta v^{\bar j}$ lies in the null space of the Ricci tensor. Then
\[0=D_k(R_{i\bar j}D_{\bar k}v^i)=D_kR_{i\bar j}D_{\bar k}v^i + R_{i\bar j}D_kD_{\bar k} v^i,\]
and
\[0=D_{\bar k}(R_{i\bar j}D_{k}v^i)=D_{\bar k}R_{i\bar j}D_{k}v^i + R_{i\bar j}D_{\bar k} D_{k} v^i\]
and so
\[D_kR_{i\bar j}D_{\bar k}v^i + D_{\bar k}R_{i\bar j}D_{k}v^i=0.\]
 Then
\[0=\Delta (R_{i\bar j}v^i)=(\Delta R_{i\bar j}) v^i +{1\over 2}(D_kR_{i\bar j}D_{\bar k}v^i +D_{\bar k}R_{i\bar j}
D_k v^i)\] and hence $(\Delta R_{i\bar j}) v^i =0$. Then
\[ 0={\partial \over {\partial t}}(R_{i\bar j}v^i)=R_{i\bar j}{\partial v^i\over {\partial t}} +(\Delta R_{i\bar j}
+R_{i\bar j k\bar l}R_{l\bar k})v^i. \] Now $R_{i\bar j k\bar l}R_{l\bar k}v^i=0$ whenever $v$ is in the null space of $R_{i\bar j}$. Thus $R_{i\bar j}v^i=0$, and ${\partial v}/{\partial t}$ lies in the null space of the Ricci tensor
as well. This shows the null space of the Ricci tensor is invariant in time. Therefore, either the null space of the Ricci tensor is trivial, or, by the De Rham decomposition theorem (see e.g. Theorem 8.1 in [12]) and induction on the dimension of the null space, the underlying complex manifold $X$ splits holomorphically isometrically into a product of some flat complex Euclidean space ${\bf C}^k$ and another factor $N^{n-k}$ whose Ricci tensor is everywhere positive. This completes the proof of Theorem 2.1.\\

\noindent {\bf Remark}: A similar argument of the proof has been used recently by Ni-Tam in their
work [14] on a Liouville type theorem for plurisubharmonic functions.

\section{Type II singularity models and K\"ahler-Ricci solitons}

In this section, we collect some known results about Type II singularity models and translating K\"ahler-Ricci solitons.

An important class of Type II singularity models is given by {\em translating K\"ahler-Ricci soliton} (KRS), which is an eternal solution $g$ moving along the K\"ahler-Ricci flow (1) under a one-parameter family of biholomorphisms of a noncompact complex manifold $X$ generated by some holomorphic vector field $V$ on $X$. That is we have $g(t)=\phi^{*}(t)g(o)$, where $\phi (t)=\exp\{-tV\}$ is the one-parameter family of automorphisms of $X$. Equivalently, for each $t$, we have
$$R_{i \bar j}=L_{V}g_{i \bar j},$$
the Lie derivative of $g$ in the direction of $V$. Thus in local
holomorphic coordinates, translating KRSs are characterized by
equations
\begin{eqnarray*}
 R_{i\bar j}=D_{\bar j} V_i=D_{i}V_{\bar j}, \qquad \mbox {and} \qquad D_jV_i=D_{\bar j}V_{\bar i}=0.
\label{3.1}
\end{eqnarray*}
Note that the condition $D_jV_i=0$ is equivalent to saying that
the vector field $V$ is holomorphic. If the vector field $V$ is the gradient of a real-valued smooth function $f$ on $X$ so that $V^{i}=g^{i\bar j}\partial_{\bar j}f$, then we call $g$ a {\it gradient translating KRS} and $f$ a {\it potential function} of the soliton. In this case, the above soliton equation becomes
\[ R_{i\bar j}=D_iD_{\bar j}f, \qquad \mbox {and} \qquad D_iD_j f=0. \]
It turns out all Type II singularity models arise this way:\\

\noindent{\bf Proposition 3.1} (Cao [3]) {\em Any Type II singularity model $(X, g)$
of the K\"ahler-Ricci flow (1) with nonnegative holomorphic bisectional curvature
and positive Ricci curvature is necessarily a translating  K\"ahler-Ricci soliton. Furthermore, if $X$ is simply connected, then $g$ is a gradient soliton.}\\

The proof of Proposition 3.1 follows from our Li-Yau-Hamilton estimate [1,3] (also previously called Harnack estimate) for the K\"ahler-Ricci flow and the strong maximum principle argument.

In complex dimension one there is only one (up to scaling) translating KRS of positive curvature, called cigar soliton, found by Hamilton [10]. The metric is defined on the complex plane ${\bf C}$ and can be written explicitly, at $t=0$, as
\[ ds^2={|dz|^2\over {1+|z|^2}}.\]
The cigar soliton has maximal curvature at the origin and is
asymptotic to a flat cylinder at the infinity.
Later in [2], the author found for each $n\geq 2$ a translating KRS on ${\bf C}^n$, invariant under the unitary group and of positive sectional curvature. Furthermore, this rotationally symmetric soliton has the following geometric properties: curvature at geodesic distance $s$ from the origin decays like $1/s$; while volume of geodesic ball of radius $s$ centered at the origin grows like $s^n$. In general, the special nature of translating KRSs also allow us to draw some very nice conclusions about both its complex analytic and geometric properties:\\

\noindent{\bf Proposition 3.2} (Cao-Hamilton [5]) {\em Let $(X^n, g)$ be a translating gradient KRS with positive Ricci curvature such that the scalar curvature
$R$ assumes its maximum $1$ in space-time. Let $f$ be a potential function of $g$. Then $f$ is a strictly convex exhaution function on $X$. In particular, $X$ is a Stein manifold diffeomorphic to $R^{2n}$.}\\

A fact important in the proof of Proposition 3.2 is that the potential function $f$ satisfies the equation
\[ R+|Df|^2=1.\]

\noindent{\bf Proposition 3.3} (Cao [4]) {\em Let $(X^n, g)$, $n\ge 2$, be a translating KRS with bounded and nonnegative bisectional curvature such that the scalar curvature $R$ assumes its maximum in space-time. For any (small) constant $\epsilon>0$, if there exists a positive constant $C_{\epsilon}$ such that the $({1+\epsilon})$-asymptotic scalar curvature ratio

\[ A_{1+\epsilon}=\lim \sup_{s\rightarrow \infty} R(x)  s^{1+\epsilon}(x) < C_{\epsilon}, \]
 then $g$ must be flat. In particular, if $g$ is non-flat, then $A_{1+\epsilon}=\infty $ for any $\epsilon>0$.}\\

Proposition 3.3 is a K\"ahler analogue of a similar result of Hamilton (Theorem 20.2 in [11]) for translating Ricci solitons of positive sectional curvature in the Riemannian case.
Note in particular the asymptotic curvature ratio $A\equiv A_2$ must be infinite for any translating
KRS of nonnegative bisectional curvature.

\section{The proof of main results}

\noindent {\bf Proof of Theorem 1.1}:  Let $(X^2, g)$ be a Type II singularity model
as in the statement of Theorem 1.1. For simplicity we assume $X^2$ is simply connected, otherwise we can replace it by its universal cover. First we claim the Ricci curvature 
of $g$ must be strictly positive at all points and all time. If the Ricci
curvature is not strictly positive at some point $x_0$ and some
time $t_0$, then by the splitting theorem 2.1, $(X^2, g)$ splits
as a quotient of a product of the flat complex plane ${\bf C}$ with a complete 
Riemann surface $\Sigma$ of positive Gaussian curvature. Note $g$
restricted to $\Sigma$ remains a Type II singularity model satisfying the
local injectivity radius estimate. But the only Type II singularity
model on a Riemann surface with positive Gaussian curvature is
the cigar soliton, which does not satisfy the local injectivity
radius estimate. A contradiction. Therefore the Ricci curvature of
$g$ is strictly positive everywhere.

Now, it follows from Proposition 3.1 that the Type II limit
solution $(X^2, g)$ is a translating gradient K\"ahler-Ricci
soliton with bounded nonnegative holomorphic bisectional curvature
and positive Ricci curvature. Then by Proposition 3.3, we know
that asymptotic scalar curvature ratio $A=A_2<\infty$. Therefore, we can 
apply Lemma 22.2 of Hamilton [11] to find a sequence of points $x_j\in X$ 
going to infinity at time $t=0$, a
sequence of radii $r_j>0$, and a sequence of positive numbers
$\delta_j\rightarrow 0$ such that

\begin{enumerate}
\item[{(a)}]  $R(x,0)\leq (1+\delta _j)R(x_j,0)$ for all $x$ in
the ball ${B}_{r_j}(x_j,0)$ of radius $r_j$ around $x_j$ at time
$t=0$;

\item[{(b)}]  $r_j^2 R(x_j,0)\rightarrow \infty $;

\item[{(c)}]  if $s_j$ is the distance of $x_j$ from some origin
$O$ at time $t=0$, then $\lambda_j=s_j/r_j\rightarrow \infty $;

\item[{(d)}]  the balls $B_{r_j}(x_j,0)$ are disjoint.

\end{enumerate}

We can then blow down $(X^2, g)$ as in [11] by taking a sequence of dilations
of $(X^2, g)$ around the sequence of points $x_j$ which we take as
our new origins $O_j$, and we shrink down instead of expanding to
make the scalar curvature $R(x_j,0)$ dilate to equal $1$ at $(O_j,
0)$. The balls ${B}_{r_j}(x_j,0)$ are then dilated to the balls
centered at the origin $O_j$ of radii $\hat
r_j=r_j^2{R}(x_j,0)\rightarrow \infty $ by property (b). Property
(a) gives good bounds on the curvature in these balls at time
$t=0$, while the same bounds for $t\leq 0$ follows from the fact
that the scalar curvature $R(x, t)$ of $g$ is pointwise increasing
in time, a consequence of our Li-Yau type estimate for $R(x,t)$
(see [1,3]). The local injectivity radius estimate (2) now becomes

\[\mbox{inj}_{B_{r_j}}\left( O_j,g\right)
\geq \frac \beta {\sqrt{R(O_j,0)}}\] since it is invariant under
dilation.

Now we have everything we need to take a limit of the dilations of
$(X^2,g)$ around the points $(x_j,0)$, dilating time like distance
squared and keeping $t=0$ as $t=0$ in the new limit, which is
denoted by $(\widehat X^2,\widehat{g})$, where $\widehat g$ is a
complete solution to the K\"ahler-Ricci flow (1) on the limiting
noncompact complex surface $\widehat X^2$ with the limiting
origin point $\widehat O$ and an ancient time interval $-\infty <t
\leq \Omega$, for some $\Omega >0$ by Shi's short time existence
result [16], and with bounded and nonnegative holomorphic
bisectional curvature such that the scalar curvature $\widehat
R(x, t)\leq 1$ everywhere for $t\leq 0$ and $\widehat R(\widehat
O,0)=1$ at the origin at time $t=0$.

The next step is to produce a Ricci-flat direction at $\widehat O$ in the blow-down 
limit $(\widehat X^2,\widehat{g})$. First we prove a proposition, which is of independent interest and will be important in showing the Ricci tensor $\widehat R_{i\bar
j}$ of $\widehat g$ must have a null vector at $\widehat O\in
\widehat X^2$. \\

\noindent {\bf Proposition 4.1} {\em Let $\Theta$ denote the
second Chern form of $(X^2, g)$ at $t=0$. Then, the second Chern
number integral $$c_2(X)=\int_X \Theta$$ of $(X^2, g)$ is well defined.
In fact, $c_2(X)\leq 1$. }\\

\noindent {\bf Proof}:  It is well-known that for any K\"ahler
surface of nonnegative holomorphic bisectional curvature, the
second Chern form, or equivalently the Gauss-Bonnet-Chern
integrand, $\Theta$ of $(X^2, g)$ is pointwise nonnegative. In
fact, $\Theta$ can be expressed explicitly as
\begin{eqnarray}
\Theta &= &{1\over 8{\pi}^2} (R^2+|Rm|^2-2 |Rc|^2)d\mu \nonumber \\
         &=&{1\over 4{\pi}^2}(R_{1\bar 1 1\bar 1}R_{2\bar 2 2\bar 2}
+2 R^2_{1\bar 1 2\bar 2}+ |R_{1\bar 2 1\bar 2}|^2+2|R_{1\bar 1
1\bar 2}|^2 + 2|R_{2\bar 2 1\bar 2}|^2)d\mu.
\end{eqnarray}
Here $d\mu$ is the volume form of
$g$ at $t=0$. Thus, in order to prove $c_2(X)$ is well defined, it suffices to
prove the integrals $c_2(U_i)=\int_{U_i} \Theta$ is uniformly
bounded for some exhaution sequence of subsets $U_i$ of $X^2$.

Since $(X^2, g)$ is a gradient KRS, by Proposition 3.2, there exists a strictly convex exhaustion potential function $f$ on $X$. Taking $U_i=\{f\le i\}$, then as is well-known, we have

\[c_2(U_i)=\int_{U_i}\Theta =\chi \left(U_i\right)+\{\mbox {boundary contribution}
\} \]
where $\chi \left(U_i\right)$ denotes the Euler characteristic of the set $U_i$.
Now $\partial U_i$ has positive definite second fundamental form
and it can be shown (see e.g. Section 4 in [8]) that the boundary term above is nonpositive
. Therefore
\[c_2(U_i)\le \chi \left(U_i\right).\]
But by Proposition 3.2, $X$ is diffeomorphic to ${\bf R}^{4}$ and 
each $U_i$ is diffeomorphic to the unit ball. Thus $\chi \left(U_i\right)=1$ for all $i$ which
implies $c_2(X)\leq 1$.  \\

Now we can prove \\

 \noindent {\bf Lemma 4.2} {\em The Ricci tensor
of $\widehat g$ must have a zero eigenvalue at the origin $\widehat O\in \widehat X^2$ at time $t=0$.}\\

\noindent{\bf Proof of Lemma 4.2}: We prove by contradiction.
Suppose Lemma 4.2 is not true. Then before taking the blow-down
limit, there must exist a positive number $\varepsilon >0$ and a
subsequence, again indexed by $j$, of $\{x_j\}$ in $X$ such that
\begin{eqnarray}
\varepsilon _{j}=\frac{\sigma_{j}}{{R}(x_{j},0)}\geq \varepsilon \
\qquad \mbox{for all} \quad j=1,2,\cdots \ ,
\end{eqnarray}
where we denote by $\sigma_{j}$ the minimum of the Ricci curvature
of $g$ at $x_j$ at $t=0$.

We are going to show this leads to a contradiction by closely examing the relation
between Chern numbers $c_2$ and $c^2_1$.

On one hand, by the local derivative estimate of Shi [16] (or
Theorem 13.1 in [11]) and properties (a) and (b), we have at $t=0$
the estimate

\begin{eqnarray}
\sup \limits_{x\in B_{r_j}(x_{j},0)}\left| \nabla Rm(x,0)\right|
^2&\leq&C R^2(x_j,0)
\left( \frac 1{r_{j}^2}+ R(x_j,0)\right) \nonumber \\
&\leq&2C R^3(x_j,0)\ ,
\end{eqnarray}
where $Rm$ is the curvature tensor of $g$ and $C>0$ is a universal
constant.

For any $x\in B_{r_j}(x_{j},0)\subset X$, we obtain from (5) and
(6) that the minimum of the Ricci curvature $\sigma_{\min}(x)$ at
$x$, satisfies
\begin{eqnarray}
\label{}
\sigma_{\min }(x)&\geq&\sigma_{j}-\sqrt{2C} R^{3/2}(x_{j},0)d_0(x,x_{j})\nonumber\\
&\geq& R(x_{j},0)\left( \varepsilon -\sqrt{2C}\cdot \sqrt{R(x_{j},0)}\cdot d_0(x,x_{j})\right)\\
&\geq&\frac \varepsilon 2 R(x_{j},0)\nonumber
\end{eqnarray}
whenever
$$d_0(x,x_{j})\leq l_j=:\frac{\varepsilon}{2\sqrt{2C}\cdot \sqrt{R(x_{j},0)}}\ . $$
Thus, from property (a) and (7), there exists $j_0>0$ such that
for any $j\geq j_0$ and $x \in B_{l_j}(x_{j},0)$, we have the
estimate
\begin{eqnarray}
\frac \varepsilon 2 R(x_{j},0) \leq  \sigma_{\min}(x)  \leq 2
R(x_{j},0)
\end{eqnarray}
in the geodesic ball $B_{l_j}(x_{j_k},0)$.

Hence

\[ {Rc_{g}}^2 (x) \geq ({\varepsilon\over 2})^2 R^2(x_{j},0) d\mu \]
for all $x\in B_{l_j}(x_{j},0)$.

Therefore, 
\begin{eqnarray}
 \sum\limits_{j=j_0}^\infty \int_{B_{l_j}(x_{j},0)} {Rc_{g}}^2 &
\geq & C(\varepsilon)\sum\limits_{j=j_0}^\infty {R}^2(x_{j},0)
\cdot C_1 \left( \frac\varepsilon {2\sqrt{2C}\cdot \sqrt{R(x_{j},0)}}\right) ^{4} \nonumber\\
& = & C(\varepsilon )\sum\limits_{j=j_0}^\infty
\frac{C_1\varepsilon ^{4}}{C^2}\\
& = & +\infty \nonumber
\end{eqnarray}

On the other hand, if we diagonalize the Ricci tensor at each
point $x$, we have
\begin{eqnarray}
{Rc_g}^2(x)&=&{1\over \pi^2} R_{1\bar 1}(x)R_{2\bar 2}(x) d\mu \nonumber\\
&=&{1\over {\pi}^2}(R_{1\bar 1 1\bar 1}R_{2\bar 2 2\bar 2}+R_{1\bar 1 1\bar 1}R_{1\bar 1 2\bar 2}+ R^2_{1\bar 1 2\bar 2}+ R_{1\bar 1 2\bar 2}
R_{2\bar 2 2\bar 2}+R^2_{2\bar 2 2\bar 2}) d\mu.
\end{eqnarray}

From (4), (5), and (10), it is easy to check that there exists a
constant $C'(\epsilon)$ such that
$$ {Rc_g}^2(x_j)\leq C'(\epsilon) \Theta (x_j)$$
for all $j$. It then follows that
$$\sum\limits_{j=j_0}^\infty \int_{B_{l_j}(x_{j},0)} {Rc_{g}}^2\leq
C'(\epsilon)\sum\limits_{j=j_0}^\infty
\int_{B_{l_j}(x_{j},0)}\Theta < C'(\epsilon)
c_2(X)<C'(\epsilon),$$ which contradicts (9). Thus the proof of Lemma 4.2 is completed.\\

Now, it follows from Lemma 4.2  and Theorem 2.1 that, since
$\widehat g$ is nonflat, the universal cover of $\widehat X^2$
splits holomorphically isometrically as a product of the flat
complex plane ${\bf C}$ and a Riemann surface $\Sigma$ of positive
Gaussian curvature. The factor $\Sigma$ may not be yet of Type I
or Type II, but we can take a further limit of dilation, as in
[11], also by shrinking, to get yet another limit $N$ which will
be either of Type I or Type II. However, $N$ cannot be Type II,
because otherwise $N$ has to be the cigar soliton, but cigar
soliton does not satisfy the local injectivity radius estimate. So
$N$ must be of Type I. Then by a result of Hamilton [11], the only
Type I ancient solutions to the Ricci flow on a Riemann surface
which are complete with bounded curvature are the (round) Riemann
sphere, or the flat complex plane ${\bf C}$ and its quotients.
Since $N$ has positive curvature, it must be the Riemann sphere
${\bf CP}^1$.

Now a sequence of dilations of $(X^2, g)$ converges to a quotient
of $(\widehat X^2, \widehat g)$, which splits into ${\bf C}\times
\Sigma$, and a sequence of dilations of $C\times \Sigma$ converges
to ${\bf C}\times CP^1$. Again, since a dilation of a dilation is a
dilation, and a limit of limits is a limit by picking an
appropriate subsequence.
Thus a limit of dilations of $(X^2, g)$ converges to ${\bf C}\times CP^1$. \hfill $\Box$ \\

\medskip

\noindent {\bf Remark}: In [7], when curvature operator is nonnegative, the authors used
the finiteness of the Euler number (or equivalently the second Chern number), to produce 
a zero holomorphic sectional curvature direction at $\widehat O$
at time $t=0$. There is an alternative proof of Theorem 1.1 along the same line. Here, once we have Proposition 4.1, we can do the same thing
to get a flat holomorphic sectional curvature direction for $\widehat g$.  But by a result of Mok [13] that is also a Ricci flat direction for $\widehat g$. 

\medskip

\noindent{\bf Proof of Corollary 1.1}: 
Suppose we have a $2$-dimensional Type II singularity model $(X^2,
g)$ with $\nu_X>0$. Using the convex exhaustion function $f$ in
Proposition 3.2 and modifying an argument of Cheeger-Gromov-Taylor
[6], it follows that the local injectivity radius estimate (ii)
holds for $(X^2, g)$. Hence, according to Theorem 1.2, a sequence
of dilations of $g$ converges to a limit $\widehat g$ on some
noncompact complex manifold $\widehat X$. Furthermore $(\widehat
X^2, \widehat g)$ is a quotient of the product ${\bf C}\times {\bf
CP}^1 $ of flat complex plane ${\bf C}$ and the complex projective
plane ${\bf CP}^1$ so that $\nu_{\widehat X}=0$. On the other
hand, we must have $\nu_{\widehat X}>0$ since the condition of
having Euclidean volume growth is clearly preserved under
dilations and taking limits.  A contradiction. \hfill $\Box$

\medskip

\noindent{\bf Proof of Theorem 1.2}: We can proceed similarly 
as in the proof of Theorem 1.1. Let $g$ be the ancient
solution to the K\"ahler-Ricci flow (1) satisfying all the
assumptions in Theorem 1.2. Since the asymptotic scalar curvature
ratio $A=\infty$, we can apply Lemma 22.2 of Hamilton [11] 
and get blow-down limit $(\widehat X^n, \widehat
g)$ as in the proof of Theorem 1.1.\\

\noindent {\bf Lemma 4.3} {\em The Ricci tensor of $\widehat g$
must have a zero eigenvalue at the origin $\widehat O\in \widehat X^n$ at time $t=0$.}\\

\noindent{\bf Proof of Lemma 4.3}: We again prove by contradiction.
Suppose Lemma 4.3 is not true so that condition (5) in the proof of Lemma 4.2
holds.  We are going to show this leads to a contradiction to assumption
(iii).

On one hand, as in the proof of Lemma 4.2, we have the
estimate
\begin{eqnarray*}
\frac \varepsilon 2 R(x_{j},0) \leq  \sigma_{\min}(x)  \leq 2
R(x_{j},0)
\end{eqnarray*}
in the geodesic ball $B_{l_j}(x_{j_k},0)$ and hence

\[ {Rc_{g}}^n (x) \geq ({\varepsilon\over 2})^n R^n(x_{j},0) d\mu \]
for all $x\in B_{l_j}(x_{j},0)$. Here $d\mu$ is the volume form of
$g$ at $t=0$. It then follows that
\begin{eqnarray*}
 c_1^n(X)&\geq & \sum\limits_{j=j_0}^\infty \int_{B_{l_j}(x_{j},0)} {Rc_{g}}^n \\
& \geq & C(\varepsilon)\sum\limits_{j=j_0}^\infty {R}^n(x_{j},0)\cdot C_1 \left( \frac\varepsilon {2\sqrt{2C}\cdot \sqrt{R(x_{j},0)}}\right) ^{2n} \\
& = & C'(\varepsilon )\sum\limits_{j=j_0}^\infty
\frac{C_1\varepsilon ^{2n}}{C^n}\\
& = & +\infty,
\end{eqnarray*}
which is a contradiction to assumption (iii), and the proof of
Lemma 4.3 is completed.\\

Now the blow-down limit $(\widehat X^n, \widehat g)$ has a
Ricci-flat direction, hence splitting theorem 2.1 applies and the
universal cover of $\widehat X$ splits holomorphically
isometrically as a product of flat ${\bf C}^l$ and a factor $N^{n-l}$ of
nonnegative bisectional curvature and positive Ricci curvature.
The factor $N$ may not be yet of Type I or Type II, but we can
take a further limit of dilation, as in [11], also by shrinking,
to get yet another limit $\tilde N$ which will be either of Type I
or Type II. If $\tilde N$ is of Type I, or Type II with
$c_1^{n-l}(\tilde N)=\infty$, then we are done. If $\tilde N$ is
of Type II but with $c_1^{n-l}(\tilde N)<\infty$, then we can
repeat the dimension reduction process above to split out more
flat factors until we arrive at some factor $\widehat N^{n-k}$
($k\ge l$) which is either of Type I, or of Type II with
$c_1^{n-k}(\widehat N)=\infty$.

Now a sequence of dilations of $(X^n, g)$ converges to a quotient of
$(\widehat X^n, \widehat g)$, which splits into ${\bf C}^{l}\times N^{n-l}$, and
a sequence of dilations of ${\bf C}^l\times N^{n-l}$ converges to ${\bf C}^l\times \tilde N^{n-l}$,
and a sequence of dilations of ${\bf C}^l\times \tilde N^{n-l}$ converges
to ${\bf C}^k\times \widehat N^{n-k}$. Since a dilation of a dilation is a dilation,
and a limit of limits is a limit by picking an appropriate subsequence.
Thus a limit of dilations of $(X^n, g)$ converges to ${\bf C}^k\times \widehat N^{n-k}$, where $\widehat N^{n-k}$ is either of Type I, or of Type II
with $c_1^{n-k}(\widehat N)=\infty$. \hfill $\Box$ \\

\noindent Institute for Pure and Applied Mathematics at UCLA, Los Angeles, CA 90095 \&

\noindent Department of Mathematics, Texas A\&M University, College Station, TX 77843 \\
{\em E-mail}:  hcao@ipam.ucla.edu  or  cao@math.tamu.edu\\

\end{document}